
\magnification = \magstep1
\hoffset 0.3 true in
\voffset 0.3 true in
\hsize = 5.3 true in
\vsize = 8.5 true in
\baselineskip =15  pt
\parskip 10 pt

\mathsurround=1pt

\hbadness 10000
\vbadness 10000


\def \acapo {\hfill\break}
\def \fine  {\hfill\vrule height5 true pt width 5 true pt depth 0 pt}
\def \sp    {\hbox{\ \ \ }}

\def \impl   {\hbox{ $\Longrightarrow$ }}
\def \se     {\hbox{ $\Longleftarrow$ }}
\def \sse    {\hbox{ $\Longleftrightarrow$ }}


\def\bin (#1 #2){{{#1}\choose{#2}}}

\def \ch#1    {\overline{ #1 }}

\def \ogni   {\ \forall~}
\def \cto    {\subseteq}
\def \stella {{^*}}
\def \div    {\big|}
\def \meno   {\backslash}
\def \vuoto{ $\O$ }

\def\b{\beta}
\def\s{\sigma}

\def\B{{\bf B}}
\def\L{{\bf L}}
\def\S{{\bf S}}
\def\X{{\bf X}}
\def\N {\bf N}

\def\gen{{\bf G}}
\def\bas{{\bf G}_k}

\def \ProclaimPar #1{\par {\parindent=-14 true pt   #1}}
\def \CasePar #1{    \par {\parskip0cm #1}}

\def \Thm#1   {{\bf\ProclaimPar Theorem \section.#1.} }
\def \Prop#1  {{\bf\ProclaimPar Proposition \section.#1.} }
\def \Lemma#1 {{\bf\ProclaimPar Lemma \section.#1.} }
\def \Cor#1   {{\bf\ProclaimPar Corollary \section.#1.} }
\def \Pf      {{\it\ProclaimPar Proof.} }
\def \Def     {{\bf\ProclaimPar Definition.} }
\def \Rem     {{\bf\ProclaimPar Remark.} }

\def\section{0}

\ \vskip1true in

\centerline{UPPER BOUNDS FOR THE BETTI NUMBERS}
\centerline{OF A GIVEN HILBERT FUNCTION}
\acapo\acapo
\centerline{\bf Anna Maria Bigatti}
\acapo
\centerline{Dipartimento di Matematica dell'Universit\`a di Genova}
\centerline{Via L.B. Alberti 4, I-16132 Genova Italy}
\centerline{E-mail Bigatti@UniMat.To.CNR.It}

\acapo\acapo

Let $R:=k[X_1,\dots,X_N]$ be the polynomial ring in $N$ indeterminates over a
field $k$ of characteristic 0 with $\deg(X_i)=1$ for $i=1,\dots,N$, and let
$I$ be a homogeneous ideal of $R$.  The Hilbert function of $I$ is the function
from $\N$ to $\N$ which associates to every natural number $d$ the dimension
of $I_d$ as a $k$-vectorspace.

$I$ has an essentially unique minimal graded free resolution
$$
0 \longrightarrow L_m
{\buildrel d_m\over\longrightarrow} L_{m-1}
{\buildrel d_{m-1}\over\longrightarrow}  \dots
{\buildrel d_2\over\longrightarrow} L_1
{\buildrel d_1\over\longrightarrow} L_0
{\buildrel d_0\over\longrightarrow} I
 \longrightarrow 0
$$
which is characterized, among the free graded
resolutions, by the condition
$$
d_q(L_q) \cto (X_1,\dots,X_N)L_{q-1}\sp \ogni
q\ge 1
$$

And therefore the Betti numbers, which are defined by
$$\b_q(I):={\rm rank}L_q$$
are invariants of $I$.

{}From Macaulay [M] (see also Robbiano [R]) it follows that a lex-segment ideal
has the greatest number of generators (the 0-th Betti number $\b_0$) among all
the homogeneous ideals with the same Hilbert function.

In this paper we prove that this fact extends to every Betti number, in the
sense that
 all the Betti numbers of a lex segment ideal are bigger than or equal to
the ones of any homogeneous ideal with the same Hilbert function.

Section 1 gives some useful notation and definitions and many simple
properties of Borel normed sets.

In Section 2 a Theorem is derived (Theorem 2.1) which is our  main
tool in comparing lex-segment and Borel normed sets.

In Section 3, using a result due to Eliahou and Kervaire [E-K], we compare
lex-segment and Borel normed ideals, and then, using some results
due to  Galligo [Ga] and  M\"oller-Mora [M-M], we compare
lex-segment and homogeneous ideals.

Section 4 gives the formula which computes the Betti numbers of the lex-segment
ideal, given its Hibert function, and these are the sharp upper bounds for the
Betti numbers of any homogeneous ideal with the same Hilbert function.

\acapo
{\bf 1.Some remarks on Borel normed sets.}

\def\section{1}

\bigskip

\ProclaimPar
{\bf Notations.} Let
 $\X_N$ denote the set of indeterminates $\{X_1,\dots,X_N\}$; then
 $(\X_N)^D$  indicates the set of all monomials of degree $D$ in $\X_N$.
\acapo
Let $\S$ be a
subset  of $(\X_N)^D$; then    $\X_N\S$ denotes the multiples of
$\S$ of degree $D+1$, i.e. $\X_N\S=\cup_{T\in \S} \{X_1T,\dots,X_NT \}$.
 \acapo
If $T=X_1^{t_1}\dots X_N^{t_N}$, then we denote by
 $m(T)$ := $\max\{i \;|\; t_i>0\}$, i.e.  the largest index of the
indeterminates actually occuring in $T$.

\Def
A set of monomials $\S\cto(\X_N)^D$ is {\bf Borel normed} if
$
T\in \S$ implies $X_i{T\over X_j}\in \S$ for all $j$ such that $X_j$
divides $T$ and for all $i<j$.

\Def
On $(\X_N)^D$ we will use the {\bf lexicographic order}, i.e.
if $T=X_1^{t_1}\dots X_N^{t_N}$ and $T'=X_1^{s_1}\dots X_N^{s_N}$ are two
monomials
in $(\X_N)^D$ we will say that $T>T'$ if $t_1=s_1,\dots,t_{i-1}=s_{i-1}$ and
$t_i>s_i$.
\acapo
Note that it is a total ordering and then there exists the minimum of every
subset of $(\X_N)^D$.

\Lemma 1
Let $\S$ be a Borel normed set.
\acapo
Then $X_i(\min \S) \in\X_N(\S\meno \{\min \S \})$   \sse $i<m(\min \S)$.

\Pf
Let $T:=\min \S$.
\CasePar
{}
\item {`\impl'}:
  If $X_iT\in\X_N(\S\meno \{T\})$ then
  $X_iT=X_jT'$ where $T'\in \S\meno \{T\}$
      \acapo
  Thus $T'>T$, hence $i<j$.
    \acapo
  On the other hand $X_j$ divides $T$
  and so $j\le m(T)$.
    \acapo
  Then $i<m(T)$.
\item {`\se'}:
  If $i<m(T)$ then,
  since $\S$ is Borel normed, $T':=X_i{T\over X_{m(T)}}\in \S \meno \{T \}$.
  Hence $X_iT=X_{m(T)}T'\in \X_N(\S\meno \{T\})$.\fine

\Prop 2
Let $\S$ be a Borel normed set.
\acapo
Then
$\X_N \S=\cup _{T\in \S}\{X_{m(T)}T,\dots, X_NT\}$ and this is a disjoint union
(i.e. $\{X_{m(T)}T,\dots, X_NT\} \cap
 \{X_{m(T')}T',\dots, X_NT'\}   =  \vuoto$ $\ogni T'\ne T$).
\Pf
    By induction on the cardinality of the set:
\CasePar {}
\item {$|\S|=1$:}
      Since $\S$ is Borel normed it follows that $\S= \{ X_1^D \} $.
      Then $\X_N \S=\{X_1X_1^D,\dots,X_NX_1^D\}$.
\item {$|\S|>1$:}
      Let the thesis be true if the cardinality is smaller then $|\S|$.
        \acapo
      Let $T':=\min \S$.
      Then
      $\X_N \S=\{X_1T',\dots,X_NT'\} \cup \X_N(\S\meno\{T'\})$.
      But from Lemma 1.1
      we have  $X_iT'\in \X_N(\S\meno\{T'\})$ $\ogni i<m(T')$.
      Hence
      $\X_N \S=\{X_{m(T')}T',\dots,X_NT' \} \cup \X_N(\S\meno\{T'\})$.
      \acapo
      Note that $\S \meno \{ T' \} $ is a Borel normed set and
      then, by the inductive hypothesis:
      $\X_N ( \S \meno \{ T' \} ) =
                \cup _{T\in \S\meno \{ T' \} }\{X_{m(T)}T,\dots, X_NT\}$.
       \acapo
      Therefore
       $\X_N \S=\cup _{ T\in \S } \{ X_{m(T)}T,\dots, X_NT \}$.
         \acapo
      Moreover from Lemma 1.1  we have that if $i\ge m(T')$ then  $X_iT'\not\in
      \X_N(\S\meno\{T'\})$
      and so
      $\{X_{m(T')}T',\dots, X_NT'\}\cap\X_N(\S\meno\{T'\})=\vuoto$.
      \acapo
      So, in particular,
      $\{X_{m(T')}T',\dots,X_NT'\}\cap\{X_{m(T)}T,\dots,X_NT\}=\vuoto$ $\ogni
      T>T'$ .
\fine

\Def Define
$$ m_i(\S) := \big| \  \{  T\in\S \;|\; m(T)=i \} \ \big|$$
i.e. the number of the elements of $\S$ which ``finish'' with $X_i$,
and similary
$$ m_{\le i}(\S) := \big| \ \{  T\in\S \;|\; m(T)\le i \} \ \big|$$
i.e. the number of the elements of $\S$ in the first $i$ indeterminates.

\Prop 3 Let $\S$ be a Borel normed set, then
\CasePar
{}
\item {i)}  $m_i(\X_N\S)=m_{\le i}(\S)$.
\item {ii)} $|\X_N\S|=\sum_{i=1}^N m_{\le i}(\S)$.

\Pf
\CasePar
{}
\item {i)} From Proposition 1.2 we have that
   $\cup _{T\in \S} \{ X_{m(T)}T,\dots, X_N T\}$
   is a disjoint union.
   Then,
   in such a representation of $\X_N \S$
   every monomial $T$  with $m(T)=i$
   can be uniquely expressed as
   $X_iT'$  where $T'\in\S$ and
   $i\ge m(T')$.
   Therefore there exists a 1-1 correspondence between the monomials $T$ in
$\X_N\S$
   with $m(T)=i$ and the monomials $T'$ in $\S$ with $m(T')\le i$.
   Hence
   $m_i(\X_N\S)=m_{\le i}(\S)$.
\item {ii)}  $|\X_N\S|=\sum_{i=1}^N m_i(\X_N\S)=\sum_{i=1}^N m_{\le i}(\S)$.
\fine

\Def A set of monomials $\S\cto(\X_N)^D$ is {\bf lex-segment} if
$$
T\in \S \hbox{ and } \ T'>T \Longrightarrow T' \in \S
$$
or equivalently $\S$ is a lex-segment set if and only if  $\S=\{ T \;|\;
T\ge\min
\S \}$.

\Rem $\S$ is a lex-segment set \impl $\S$ is a Borel normed set.

\Def We can uniquely decompose $\S$, with respect
to $X_N$, as follows
$$
\S= \ \S_0 \ \cup \ X_N\S_1 \ \cup \ X_N^2\S_2 \ \cup \dots  \cup \ X_N^D\S_D
\ $$
where the $\S_d$'s are sets of monomials in $N-1$ indeterminates.
\acapo
More precisely  $$\S_d\cto (\X_{N-1})^{D-d}$$

\Rem It is easy to see that
 $\ogni i<N \sp m_{\le i}(\S)=m_{\le i}(\S_0)$.
\acapo
In particular, $m_{\le N-1}(\S)=|\S_0|$.

\Prop 4
\CasePar{}
\item {i)}
Let $\S$ be a Borel normed (lex-segment) set in $N$ indeterminates. For every
$d$ \  $\S_d$ is a Borel normed (lex-segment) set in $N-1$ indeterminates.
\item {ii)}
If $\S$ is a Borel normed (lex-segment) set in $(\X_N)^D$ then
  $\X_N \S$ is a Borel normed (lex-segment) set in $(\X_N)^{D+1}$.

\Pf Easy exercise.

\Lemma 5  $\S$ is a Borel normed set \sse
   $\S_d$ is a Borel normed set $\ogni d$, \ and \
$\X_{N-1} \S_d\cto \S_{d-1}$  $\ogni d>0$.

\Pf
\CasePar{}
\item{`\impl':}
   $\S_d$ is Borel normed set follows from Proposition 1.4.i.
      \acapo
   Then it remains to prove that $\X_{N-1} \S_d\cto \S_{d-1}$  $\ogni d>0$.
      \acapo
   Let $T\in \S_d$ i.e. $X_N^dT\in \S$.
   Since $\S$ is Borel normed we have:
   $ X_N^{d-1} X_i T = X_i{X_N^d T\over X_N}\in \S$  $\ogni i<N$.
   \acapo
   Hence
   $X_i T \in \S_{d-1}$  $\ogni i<N$.
   Thus $\X_{N-1} \S_d\cto \S_{d-1}$.
\item{`\se':}
   We need to prove:
   $T\in\S$ \impl $X_i{T\over X_j}\in \S$ $\ogni i<j$ and $X_j\div T$:
      \acapo
   Let $T = X_N^d T'$ where $T'\in \S_d$. Then:
 \acapo
    Since $\S_d$ is a Borel normed set
    it follows that
      \acapo
    $X_i{ T  \over X_j} = X_N^d \Big( X_i{ T' \over X_j}\Big)  \in \S$
                                       \   $\ogni i<j<N$, $X_j\div T$.
 \acapo
    Since $\X_{N-1} \S_d\cto \S_{d-1}$ $\ogni d>0$ (and then $X_N\div T$)
    it follows that
      \acapo
    $X_i{ T \over X_N} = X_N^{d-1}( X_i T')\ \in \S$ \  $\ogni i<N$.
\fine

\Def Given a set $\S\cto (\X_N)^D$ we can uniquely define a
corresponding {\bf lex-segment with respect to $X_N$} (denoted
 $\S\stella$) as follows:
    \acapo
Recall that $\S_d\cto(\X_{N-1})^{D-d}$ then denote by $\S\stella_d$
the lex-segment set in $(\X_{N-1})^{D-d}$ with $|\S\stella_d|=|\S_d|$,
i.e. the set of the greatest $|\S_d|$ monomials in $(\X_{N-1})^{D-d}$.
    \acapo
Then $\S\stella$ :=
$ \ X_N^0 \S\stella_0 \ \cup \dots \cup \ X_N^D \S\stella_D $.

\Rem $m_{\le N-1}(\S)=|\S_0|=|\S\stella_0|=m_{\le N-1}(\S\stella)$.

\Lemma 6
Let $\S$ be a Borel normed set.
If   $m_{\le i}(\S\stella_d)\le m_{\le i}(\S_d)$ $\ogni i\le N-1$ and $\ogni
d$,
then $\S\stella$ is a Borel normed set.

\Pf
By Lemma 1.5 it sufficies to show that
 $\S\stella _d$ is a Borel normed set $\ogni d$
  \  and \
 $\X_{N-1}\S\stella _d\cto \S\stella _{d-1}$ $\ogni d>0$.
    \acapo
The fact that $\S\stella_d$ is a Borel normed set is obvious since
$\S\stella_d$
is a lex-segment set.
    \acapo
It remains to prove $\X_{N-1}\S\stella _d\cto \S\stella _{d-1}$:
  \acapo
{}From Proposition 1.3.ii it follows  that
for every Borel normed set $\S$ \
$|\X_N\S|=\sum_{i=1}^N m_{\le i}(\S)$.
But, by hypothesis we have
\acapo
$|\X_{N-1}\S\stella _d| =
 \sum_{i=1}^{N-1} m_{\le i}(\S\stella_d)\le
 \sum_{i=1}^{N-1} m_{\le i}(\S _d) =
 |\X_{N-1}\S_d|$ $\ogni d$.
Since $\S$ is Borel normed,
it follows from Lemma 1.5  $\X_{N-1}\S_d  \cto  \S_{d-1}$ $\ogni d>0$.
Hence
$$|\X_{N-1}\S\stella_d|\le|\X_{N-1}\S_d|\le|\S_{d-1}|=|\S\stella_{d-1}|
\ \ogni d$$
Thus, since $\X_{N-1}\S\stella_d$ and $\S\stella_{d-1}$ are lex-segments sets
(Proposition 1.4), we have
$$\X_{N-1}\S\stella_d\cto \S\stella_{d-1} \ \ogni d$$
\fine

\Rem We will see (in  Theorem 2.1)    that  the hypothesis of this Lemma
are always verified. Thus for every  Borel normed set $\S$ we have that
$\S\stella$ is Borel normed.

\Def Let $T=(X_1^{t_1},\dots, X_N^{t_N})\in (\X_N)^D$. Define the
{\bf corresponding monomial} $\ch T $ in $(\X_{N-1})^D$ as follows:
$$\ch T :=(X_1^{\ch t_1 },\dots, X_{N-1}^{\ch t_{N-1} })$$
$$\hbox{ where \sp }
\ch t_i :=t_i  \sp \ogni i<N-1 \hbox{ \sp  and \sp  } \ch t_{N-1}
:=t_{N-1}+t_N$$
or equivalently
$$\ch T := T \left( { X_{N-1} \over X_N } \right)^{t_N} $$

\Lemma 7
\CasePar
{}
\item {i)}
  If $T,T'\in (\X_N)^D$ \ $T\le T'$
  then $\ch T \le \ch T' $.
\item {ii)}
  Let $\S$ be a Borel normed set
  then \ $\ch {\min \S} = \min \S_0$.

\Pf
\CasePar
{}
\item {i)}
  Let  $T=(X_1^{t_1},\dots, X_N^{t_N})$, $T'=(X_1^{s_1},\dots,X_N^{s_N})$ and
  $T<T'$, then we have $t_i=s_i$ $\ogni i<j$ \ and \  $t_j<s_j$.
      \acapo
  Note that $j\ne N$ since otherwise $D=\deg T = \sum t_i <\sum s_i =\deg
T'=D$.
  \acapo
  Then:
  $\ch t_i =t_i=s_i=\ch s_i $  $\ogni i< j$ \  and, (two cases)
    \acapo
  if $j=N-1$:
     $\ch t_{N-1}  =\deg T-\sum_{i=1} ^{N-2} t_i=\ch s_{N-1} $ \
     then  $\ch T =\ch T' $;
         \acapo
  if $j<N-1$:
     $\ch t_j =t_j < s_j = \ch s_j $ \ then  $\ch T <\ch T' $.
\item {ii)}
  Obviously $\min \S\le \min \S_0$.
  \acapo
  Hence from i) it follows that $\ch {\min \S}
  \le \ch {\min \S_0} = \min \S_0$.
    \acapo
  On the other hand, since $\S$ is Borel normed, $\ch {\min\S} \in \S_0$.
  \acapo
  Thus $\ch {\min\S} = \min \S_0$.
 \fine

\acapo
{\bf 2. Comparisons between lex-segment  and Borel normed sets.}

\def\section {2}

\bigskip

\Thm 1 Let  $\L$ be a lex-segment set and $\B$  a Borel normed set in
$(\X_N)^D$ such that $|\L|\le|\B|$. Then
$$m_{\le i}(\L)\le m_{\le i}(\B) \sp \sp \ogni i=1,\dots,N$$

 \Pf
By induction on the number of indeterminates:
\CasePar {}
\item {$N=2$:} $m_{\le 2}(\L) =|\L|\le|\B|= m_{\le 2}(\B)$ \sp and \sp
               $m_{\le 1}(\L) =     1     = m_{\le 1}(\B)$.
\item {$N>2$:}
Inductive hypothesis: let the thesis be true in $N-1$
indeterminates, and then study for all $i=1\dots N$  the relations between
$m_{\le i}(\L)$  and $m_{\le i}(\B)$.
\itemitem { $i=N$:}
      $m_{\le N}(\L)=|\L|\le|\B|=m_{\le N}(\B)$.
\itemitem {$i=N-1$:}
      We need to prove:
               $$m_{\le N-1}(\L)\le m_{\le N-1}(\B)$$
      i.e.
               $$|\L_0|\le |\B_0|$$
      From the definition of the lex-segment with respect to $X_N$ we have
      $|\B_0|=|\B\stella_0|$. So it will be enough to prove
               $$|\L_0|\le |\B\stella_0|.$$
\acapo {}
     Now    $\B\stella_d$ and $\B_d$ are, $\ogni d$, respectively
  lex-segment and Borel normed sets of monomials in $N-1$ indeterminates with
the
same cardinality. Then by the inductive hypothesis it follows
$$m_{\le i}(\B\stella_d)\le m_{\le i}(\B_d) \sp
\ogni i=1,\dots , N-1 \ \ogni d$$
Then, since $\B$ is Borel normed, it follows from Lemma 1.6  that
$\B\stella$ is Borel normed.
         \acapo
      Now recall the definition of corresponding monomial in $(\X_{N-1})^D$
      and
      consider
      $\ch {\min\B\stella} $ and $ \ch {\min\L} $.
           \acapo
      Note that
      $\min \L \ge \min \B\stella $ (otherwise, since $\L$ is lex-segment,
       $\B\stella \subset \L$ and then
       $|\B|=|\B\stella|<|\L|$), therefore from Lemma 1.7.i
                 $$   \ch {\min\L} \ge \ch {\min\B\stella}   $$
      It follows from Lemma 1.7.ii  that,
      since   $\L$ and $\B\stella$ are Borel normed
      $$\ch {\min\L}          =  \min \L_0  \sp  , \sp
        \ch {\min\B\stella}   =  \min \B\stella_0       $$
      and then
      $$  \min \L_0 \ge \min \B\stella_0 $$
      Moreover, since $\B\stella$ is lex-segment  w.r.to $X_N$, we have that
      $\B\stella_0$ is a  lex-segment set in $(\X_{N-1})^D$.
      From these facts it follows that
      $$ \L_0 \cto \B\stella_0$$
      Hence
      $$ |\L_0|  \le |\B\stella_0| $$
\itemitem  {$i<N-1$:}
      From the case $i=N-1$ we have   $|\L_0|\le|\B_0|$
      where $\L_0$ and $\B_0$ are respectively lex-segment
      and Borel normed sets in $N-1$ indeterminates.
      By the inductive hypothesis
      $$m_{\le i}(\L_0)\le m_{\le i}(\B_0)
         \ \ \ogni i=1,\dots,N-1$$
      So
$$m_{\le i}(\L) = m_{\le i}(\L_0) \le m_{\le i}(\B_0) = m_{\le i}(\B)
 \ \ \ogni i<N-1$$
\fine

\Cor 2  $|\L|=|\B|$ \impl $|\X_N \L|\le|\X_N \B|$.

\Pf
$|\X_N \L| = \sum_{i=1}^N m_{\le i}(\L) \le
                   \sum_{i=1}^N m_{\le i}(\B) = |\X_N \B|$.\fine

\Def Let $\S$ be any set of monomials, then define
$$ b_q(\S) := \sum_{T\in \S}\bin (m(T)-1 q)$$

\Prop 3 Let $\B\cto(\X_N)^D$ be a Borel normed set, then
$$b_q(\B)=\bin (N-1 q)|\B|-\sum_{i=1}^{N-1}
\left[ m_{\le i}(\B)\bin (i-1 q-1) \right]$$

\Pf
$$
\sum_{T\in \B}\bin (m(T)-1 q) \ = \
\sum_{i=1}^N     \left[m_i(\B)\bin (i-1 q)\right]=
$$ $$
=\sum_{i=1}^N
\left[\Big(m_{\le i}(\B)-m_{\le i-1}(\B)\Big)\bin (i-1 q)\right] =
$$ $$
=\sum_{i=1}^N        \left[m_{\le i}(\B) \bin (i-1 q)\right]  -
 \sum_{i=0}^{N-1}    \left[m_{\le i}(\B) \bin (i   q)\right] =
$$ $$
=\bin (N-1 q) m_{\le N}(\B) +
\sum_{i=1}^{N-1}
\left[m_{\le i}(\B) \left(\bin (i-1 q) - \bin (i q)\right)\right] =
$$ $$
=\bin (N-1 q)|\B|-
\sum_{i=1}^{N-1}    \left[m_{\le i}(\B)\bin (i-1 q-1)\right]
$$
\fine

\Cor 4 Let  $\L$ be a lex-segment set and $\B$  a Borel normed set in
$(\X_N)^D$ such that \ $|\L|=|\B|$, then:
\CasePar {}
\item {i)} $b_q(\L) \ge b_q(\B)$;
\item {ii)} $b_q(\X_N \L) \le b_q(\X_N \B)$.

\Pf
{}From Theorem 2.1 we have $m_{\le i}(\L)\le m_{\le i}(\B)$ $\ogni
i=1,\dots,N$.
 \CasePar {Then:}
\item {i)}
   $$b_q(\L)=    \bin (N-1 q)|\L|-\sum_{i=1}^{N-1}
   \left[m_{\le i}(\L)\bin(i-1 q-1)\right]\ge$$
   $$\ge\bin (N-1 q)|\B|-\sum_{i=1}^{N-1}\left[m_{\le i}(\B)
   \bin(i-1 q-1)\right] = b_q(\B)$$
\item{ii)} Recall  from Proposition 1.3.i that if $\S$ is a Borel normed set
   then $m_i(\X_N\S)=m_{\le i}(\S)$. Thus:
      \acapo
   $$b_q(\X_N \L)=
   \sum_{T\in \X_N\L} \bin(m(T)-1 q)=
   \sum_{i=1}^N m_{i}(\X_N\L)\bin(i-1 q)=$$
   $$=\sum_{i=1}^N m_{\le i}(\L)\bin(i-1 q)
   \le
   \sum_{i=1}^N m_{\le i}(\B)\bin(i-1 q)=
   b_q(\X_N \B)$$
\fine

\acapo
{\bf 3. Comparisons between lex-segment and homogeneous ideals.}

\def\section {3}

\bigskip

\Def Let $I$ be  a monomial ideal in $k[X_1,\dots,X_N]$. Then we denote by
$\gen(I)$ the minimal system of generators of $I$, i.e. the set of all
monomials in $I$ which are not proper multiples of any monomial in $I$,   and
by
$\bas(I_d)$ the basis of $I_d$ as a $k$-vectorial space.

\Def A monomial ideal $I$ in $k[X_1,\dots,X_N]$ is called:
\acapo
i) {\bf Lex-segment} if $\bas(I_d)$ is a lex-segment set $\ogni d$;
\acapo
ii) {\bf Borel normed} if $T\in I$
  \impl $X_i{T\over X_j}\in I$  $\ogni i<j$  such that $X_j\div T$,
  or equivalently if  $\bas(I_d)$ is a Borel normed set $\ogni d$;
\acapo
iii) {\bf Stable} if $T\in I$ \impl $X_i{T\over X_{M(T)}}\in I$.

\Rem  If $I$ is a monomial ideal. Then
\acapo
$I$ is lex-segment \impl $I$ is Borel normed \impl $I$ is stable.

\Thm 1.Eliahou-Kervaire(1987) Let  $I$  be a  stable ideal, then
 $$\b_q(I)=\sum_{T\in \gen(I)}\bin (m(T)-1 q)$$

\Cor 2 Let $I$ be a stable ideal.
\acapo
Then
$\b_q(I) = \sum_{d>0} \big[ b_q(\bas(I_d)) - b_q(\X_N \bas(I_{d-1})) \big]$.

\Pf
{}From Theorem 3.1 we have
$\b_q(I)=b_q(\gen(I))=\sum_{d>0}b_q\Big(\big(\gen(I)\big)_d\Big)$.
Then, since
$\big( \gen(I) \big)_d =\bas(I_d)\meno \{ X_N \bas(I_{d-1}) \}$,
the thesis follows.
\fine

\Cor 3 Let  $I^\L$ be a lex-segment ideal and $I^\B$  a Borel normed
ideal with the same Hilbert function, then
$$\b_q(I^\L)\ge \b_q(I^\B)$$

\Pf
Note that for all $d$
$\bas(I^\L_d)$ and $\X_N \bas(I^\L_d)$ are lex-segment sets,
and
$\bas(I^\B_d)$ and $\X_N \bas(I^\B_d)$ are Borel normed sets.
    \acapo
Moreover  $I^\L$ and $I^\B$ have the same Hilbert function, i.e. $\ogni d$
$$|\bas(I^\L_d)|=H_{I^\L}(d)=H_{I^\B}(d)=|\bas(I^\B_d)|$$
{}From Corollary 2.4, we then have
$$b_q(\bas(I^\L_d)) \ge b_q(\bas(I^\B_d))
\hbox{\  and \ }
b_q(\X_N \bas(I^\L_d)) \le b_q(\X_N \bas(I^\B_d)) \ \ogni d$$
and, from Corollary 3.2
$$\b_q(  \gen(I^\L)  )   \ =\
\sum_{d>0}  \left[ b_q\left( \bas(I^\L_d)) \right) -
                   b_q\left( \bas(I^\L_{d-1})) \right)\right]
\ge $$
$$ \ge
\sum_{d>0}  \left[ b_q\left( \bas(I^\B_d)) \right) -
                   b_q\left( \bas(I^\B_{d-1})) \right)\right]=
 \b_q(  \gen(I^\B)  )$$
\fine

\Rem
Note that for every Borel normed ideal $I$ there exists a lex
segment ideal with the same Hilbert function as that of $I$.
In fact
let $\S_d$ be the lex-segment set in $(\X_N)^d$ with $|\S_d|=H_I(d)$.
Then,
from Corollary 2.2,
$$|\X_N\S_d| \le |\X_N \bas(I^\B)|$$
Since $I^\B$ is an ideal we have
$\X_N\bas(I^\B)_d\cto\bas(I^\B)_{d+1}$.
Thus
$$|\X_N\S_d|\le|\X_N\bas(I^\B)_d| \le H_I(d+1)=|\S_{d+1}|$$
Since $\X_N\S_d$ and $\S_{d+1}$ are lex segments we get
$\X_N\S_d\cto\S_{d+1}$.\acapo
Hence we can consider the $\S_d$'s as the basis of the part in degree $d$ of an
ideal that  is lex-segment and has the same Hilbert function of $I$.

\Thm 4.Galligo(1974)
\acapo
Let $I$ be a homogeneous ideal in
$k[X_1,\dots,X_N]$ and let  $\s$ a term-ordering. There exists a Zariski
open subset  $U\cto GL(N)$ such that for every $g\in U$, $Lt_\s(g(I))$ is
invariant under the action of the Borel subgroup $B(N)$ of $GL(N)$. In
particular, if $char(k)=0$, then $Lt_\s(g(I))$ is Borel normed.

\Rem In this way we can obtain for every homogeneous ideal $I$,  an ideal
$I^\B$
with the same Hilbert function and  the same Betti numbers as those of $I$, and
such that $Lt_\s(I^\B)$  is  Borel normed.

\Thm 5.Macaulay(1927)
\acapo
Let $I$ be a homogeneous ideal in
$k[X_1,\dots,X_N]$ and let  $\s$ a term-ordering. Then
$$H_{I}=H_{Lt_\s(I)}$$

\Rem Let $I$ be a homogeneous ideal in
$k[X_1,\dots,X_N]$. Then there exists a lex
segment ideal with the same Hilbert function as that  of $I$.
\acapo
In fact let $I^\B$ be the  ideal obtained from $I$ by a generic
change of coordinates (Theorem 3.4). We have that $Lt_\s (I^\B)$ is a Borel
normed ideal  and hence there exists a lex
segment ideal with   Hilbert function $H_{Lt_\s(I^B)}=H_{I^B}$
(Theorem 3.5).

\Thm 6.M\"oller-Mora(1983)
\acapo
Let $I$ be a homogeneous ideal in
$k[X_1,\dots,X_N]$ and let  $\s$ a term-ordering. Then
$$\b_{I}\le\b_{Lt_\s(I)}$$

\Thm 7 Let $I$ be a homogeneous ideal and let $I^\L$ be the lex-segment ideal
with the same Hilbert function as that of $I$.
Then for all $q$
                     $$\b_q(I^\L)\ge\b_q(I)$$

\Pf Let $I^\B$ be the ideal obtained by Theorem 3.4. Then
$$
    H_{I^\L} = H_I = H_{I^B}
     \hbox{ \sp and \sp }
  \b_q \big( I^\B \big)=\b_q(I)
$$
{}From Macaulay's Theorem it follows that
   $$ H_{I^B} = H_{Lt_\s( I^\B )}$$
Then, from Corollary 3.3
   $$\b_q(I^\L)\ge \b_q \big(Lt_\s(I^\B)\big)$$
{}From M\"oller-Mora's Theorem
   $$\b_q \big(Lt_\s(I^\B)\big) \ge  \b_q \big( I^\B \big)$$
and then
   $$\b_q(I^\L)\ge\b_q(I)$$
\fine

\acapo
{\bf 4. Upper Bounds for Betti Numbers.}

\def\section {4}

\bigskip

\Thm 1
Let $I$ be a Borel normed ideal and, with abuse of notation, let $I_d$ denote
$\bas(I_d)$. Then
$$\b_q(I)  =$$
$$=
\bin (N-1 q)|I_D|-\sum_{i=1}^{N-1} m_{\le i}(I_D)\bin (i-1 q-1)
- \sum_{d=1}^{D-1} \left[ \
 \sum_{i=1}^{N-1}  \left[m_{\le i}(I_d) \bin(i q) \right]\ \right]
$$
Where  $D$ is the  largest degree of a generator of $I$.

\Pf
{}From Corollary 3.2  it follows that
$$
\b_q(I)=
\sum_{d=1}^D  \left[b_q( I_d ) - b_q( \X_N I_{d-1} )\right]
=
$$
$$
=
b_q( I_D )
+ \sum_{d=1}^{D-1} \left[ b_q(    I_d   )\right]
- \sum_{d=0}^{D-1} \left[ b_q( \X_N I_d )\right]
=
$$
$$
=
b_q( I_D )
+ \sum_{d=1}^{D-1} \left[ b_q(I_d)-b_q( \X_N I_d)\right]=
$$
$$
=
b_q( I_D )
+ \sum_{d=1}^{D-1} \left[ \
 \sum_{i=1}^N  \left[m_i(   I_d  )  \bin(i-1 q) \right]
-\sum_{i=1}^N  \left[m_i(\X_N I_d)  \bin(i-1 q) \right] \  \right]
$$

     Since $I_d$ is a Borel normed set  it follows  from Proposition 1.3.i
     that $m_i(\X_N I_d)=m_{\le i}(I_d)$, $\ogni d$.
     Then $\b_q(I)  =$
$$
=
b_q( I_D )
+ \sum_{d=1}^{D-1} \left[ \
 \sum_{i=1}^N  \left[( m_i(I_d) - m_{\le i}(I_d)) \bin(i-1 q) \right]
\ \right]=
$$
$$
=
b_q( I_D )
- \sum_{d=1}^{D-1} \left[ \
 \sum_{i=1}^N  \left[m_{\le i-1}(I_d) \bin(i-1 q) \right]\ \right]=
$$
$$
=
b_q( I_D )
- \sum_{d=1}^{D-1} \left[ \
 \sum_{i=1}^{N-1}  \left[m_{\le i}(I_d) \bin(i q) \right]\ \right]=
$$
$$
=
\bin (N-1 q)|I_D|-\sum_{i=1}^{N-1} m_{\le i}(I_D)\bin (i-1 q-1)
- \sum_{d=1}^{D-1} \left[ \
 \sum_{i=1}^{N-1}  \left[m_{\le i}(I_d) \bin(i q) \right]\ \right]
$$
\fine

\Def It is well known (see Robbiano [R]) that, if $h$ and $n$ are positive
integers, then $h$ can be written uniquely in the form
$$h= \bin(h(n) n) + \bin (h(n-1) n-1) + \dots + \bin (h(i) i)$$
where
$h(n) >h(n-1) > \dots > n(i) \ge i \ge 1$.
\acapo
This unique expression is called {\bf binomial expansion}
 of $h$ in base $n$ and it is denoted by $h_n$, and define
$$\big( h_n \big) ^s_t :=
\bin(h(n)+s n+t) + \bin (h(n-1)+s n-1+t) + \dots + \bin (h(i)+s i+t)$$

\medskip

The particular  significance of the binomial
expansion of the values of the Hilbert function becomes apparent when we attend
to write an explicit formula which  computes the Betti numbers of a
lex-segment ideal:
\acapo
Let $\S$ be a lex-segment set in $(\X_N)^D$ and
let $d$ be the largest integer such that
$X_1^{D-d} X_N^d\in\S$.
\acapo
Since $\S$  is a lex-segment set, $\S$ contains all the monomials
$$X_1^{D-d} \{ X_1,\dots,X_N \} ^d$$
The number of these elements is $\bin({N+d-1} {N-1})$ which is exactly
the first binomial in the binomial expansion of $H(D)$ in base $N-1$.
\acapo
The set of the remaining monomials of $\S$ is strictly contained in
$$X_1^{D-d-1} \{ X_2,\dots,X_N \} ^{d+1}$$
Thus, we can think of it as a lex-segment set (strictly contained) in
$\{ X_2,\dots,X_N \} ^{d+1}$.
So, repeating the reasoning, we obtain the
whole  binomial expansion.
\fine

\Prop 2.(Macaulay) Let  $I$ be a  lex-segment ideal. Then
$$|\X_N \bas (I_D)|=(H_I(D)_{N-1})^1$$

\Pf
As we saw before,
the first binomial of the binomial expansion of $H(D)$ in base $N-1$, $\bin
(N+d-1 N-1)$, represents the number of monomials  in $\{ X_1,\dots, X_N \}^D$.
Thus the multiples of these elements are  a set with
$\bin (N+{(d+1)}-1 N-1)$ elements.
\acapo
And so on.
\fine

\Prop 3 Let  $I$ be a  lex-segment ideal. Then
 $$m_{\le i}(I_d)=\Big( H(d)_{N-1} \Big)^{-(N-i)}_{-(N-i)}$$
(Where $\bin(h n):=0$ if $n<0$).

\Pf
As before,
$\bin (N+d-1 N-1)$ is the number of monomials
in $\{ X_1,\dots, X_N \}^D$.
Among these, the elements which use only  the first $i$ indeterminates number
$\bin({i+d-1} {i-1})$ \ i.e.
$$\bin({N+d-1-(N-i)} {N-1-(N-i)})$$
And so on.
\fine

\medskip

\Rem If $I$ is a homogeneous ideal we can calculate the largest degree of a
generator of the  lex-segment ideal  with the same   Hilbert function. In fact,
Green [Gr] proved that  $D+1$ is the smallest integer greater then the
maximum degree of a  generator of $I$ for which $H_I(D)^1= H_I(D+1)$.
\acapo
Hence Theorem 4.1 and Proposition 4.3 give a formula which computes the Betti
numbers of a lex-segment ideal.
\acapo
They are  sharp upper bounds for homogeneous
ideals with the same Hilbert function.
\acapo
In particular, to count the first syzigies, it is possible to give a simpler
formula.

\Cor 4
$\b_1(I)=$
$$(N-1)H(D)-\big( H(D)_{N-1}\big)_{-1} +
  \sum_{d=1}^{D-1}
  \bigg[(N-1)(H(d)_{N-1})_{-1} -(H(d)_{N-1})_{-2} \bigg]$$

\Pf
$b_1(\S)$
=
$\sum_{i=1}^Nm_ i(\S) \big( i-1 \big)$
=
$\sum_{i=1}^N \big( m_{\le i}(\S) - m_{\le i-1}(\S) \big) \big( i-1 \big)$
=
$\sum_{i=1}^N \bigg[ m_{\le i}(\S) \big( i-1 \big) \bigg]
   - \sum_{i=1}^{N-1} \bigg[ m_{\le i}(\S) i \bigg]$
=
$(N-1)|\S| - \sum_{i=1}^{N-1}  m_{\le i}(\S)$.
\acapo
{}From Proposition 4.3 it is easy to see that
$\sum_{i=1}^{N-1}  m_{\le i}(\S)=(|\S|_{N-1})_{-1}$.
Then
$b_1(\S)=(N-1)|\S| -(|\S|_{N-1})_{-1}$.
\acapo
{}From Proposition 4.2 it follows that
$|\X_N\S|=(|\S|_{N-1})^1$ and from Corollary 3.2 that
$\b_1(I) = \sum_{d>0} \big[ b_1(\bas(I_d)) - b_1(\X_N \bas(I_{d-1})) \big]$.
\acapo
Hence
$$\b_1=$$
$$=\sum_{d=1}^D
[(N-1)|I_d|-(|I_d|_{N-1})_{-1}-((N-1)|\X_NI_{d-1}|-(|\X_NI_{d-1}|_{N-1})_{-1})]
=$$
$$=\sum_{d=1}^D
[(N-1)H(d)-(H(d)_{N-1})_{-1}-$$
$$((N-1)(H(d-1)_{N-1})^1-(H(d-1)_{N-1})^1_{-1})]$$
The thesis follows easily.
\fine

\acapo

\centerline{ACKNOWLEGMENTS}

\medskip

Sincere thanks go to Prof. L.Robbiano and Prof. T.Mora for their useful
suggestions.

\eject

\acapo

\centerline{REFERENCES}

\bigskip

\CasePar {}
\item{[B-C-R]}
  A.M. Bigatti, M. Caboara, L. Robbiano:
  On the Computation of Hilbert-Poincar\'e Series.
  AAECC {\bf 2} (1991) 21-33.
\item{[E-K]}
  S. Eliahou, M. Kervaire:
  Minimal Resolution of Some Monomial Ideals.
  J. Algebra {\bf 129} (1990) 1-25
\item{[Ga]}
  A. Galligo:
  A propos du th\'eor\`eme de pr\'eparation de Weierstrass,
  Functions de Plusieurs Variables Complexes.
  Lecture Notes in Mathemathics
  {\bf 409}, Berlin, Heidelberg, New York: Springer (1974),  543-579.
\item{[Gr]}
  M. Green:
  Restriction of linear series  to hyperplanes, and some results of Macaulay
  and Gotzmann, Algebraic Curves and Projective Geometry Proceedings, Trento
  (1988) Springer Lecture Notes in Mathematics {\bf 1389}
\item{[M]}
  F.S. Macaulay:
  Some properties of enumeration in the theory of modular system.
  Proc. London Math. Soc. {\bf 26} (1927), 531-555.
\item{[M-M]}
  H.M. M\"oller, T. Mora:
  New Constructive methods in classical ideal theory.
  J. Algebra {\bf 100}  (1986), 138-178.
\item{[R]}
  L. Robbiano:
  Introduction to the Theory of Hilbert Function. Queen's Papers in Pure and
  Applied Mathematics {85} (1990) B1-B26.

\end